\documentclass{amsart}
\usepackage{amssymb}
\usepackage{amsfonts}
\usepackage{latexsym}

\newtheorem{theorem}{Theorem}[subsection]
\newtheorem{lemma}[theorem]{Lemma}
\newtheorem{proposition}[theorem]{Proposition}
\newtheorem{corollary}[theorem]{Corollary} 
\theoremstyle{definition}  
\newtheorem{definition}[theorem]{Definition}

\newtheorem{remark}[theorem]{Remark}
\newtheorem{note}[theorem]{Note}

\newcommand{\Tr}{\text{Tr}}
\newcommand{\id}{\text{id}} 
\newcommand{\FPdim}{\text{FPdim}} 
\newcommand{\C}{\mathcal{C}}

\newcommand{\algint}{\text{alg.\ int.}\,}

\newcommand{\End}{\text{End}} 

\newcommand{\Hom}{\text{Hom}}

\newcommand{\Rep}{\text{Rep}}

\newcommand{\eps}{\varepsilon}

\newcommand{\I}{\mathcal{I}}

\renewcommand{\C}{\mathcal{C}}
\newcommand{\M}{\mathcal{M}}

\newcommand{\actl}{\rightharpoonup}
\newcommand{\actr}{\leftharpoonup}

\newcommand{\la}{\langle\,} 
\newcommand{\ra}{\,\rangle}

\newcommand{\1}{_{(1)}} 
\newcommand{\2}{_{(2)}} 
\newcommand{\3}{_{(3)}} 
 
\renewcommand{\I}{^{(1)}} 
\newcommand{\II}{^{(2)}}   

\begin{document}

\title{Semisimple weak Hopf algebras}

\author{Dmitri Nikshych}
\address{Department of Mathematics and Statistics,
University of New Hampshire,  Durham, NH 03824, USA}
\email{nikshych@math.unh.edu}
\date{September 18}
\begin{abstract}
We develop the theory of semisimple weak Hopf
algebras and obtain analogues of a number of classical
results for ordinary semisimple Hopf algebras.
We prove a criterion for semisimplicity
and analyze the square of the antipode $S^2$ of a semisimple
weak Hopf algebra $A$. We explain how the Frobenius-Perron
dimensions of irreducible $A$-modules and eigenvalues of $S^2$
can be computed using  the inclusion matrix associated to $A$. 
A trace formula of Larson and 
Radford is extended to a relation between the categorical 
and Frobenius-Perron dimensions of $A$. Finally, an analogue of the  
Class Equation of Kac and Zhu is established and properties of 
$A$-module algebras and their dimensions are studied.
\end{abstract} 
\maketitle  


\begin{section}
{Introduction}
In \cite{BSz1} G.~B\"{o}hm and  K.~Szlach\'anyi introduced weak Hopf algebras
as a generalization of ordinary Hopf algebras. It was observed
in \cite{NV2} that these objects also extend the theory of groupoid algebras. A
general theory of weak Hopf  algebras was subsequently developed in \cite{BNSz}.  
Briefly, a weak Hopf algebra is a vector space with structures of an algebra 
and coalgebra related to each other in a certain self-dual way.
Weak Hopf algebras naturally appear in different areas of mathematics,
including algebra, functional analysis, and representation theory
\cite{KN, NV1, EN, NTV}; we refer the reader to \cite{NV2} for a survey 
of the subject. 

In this paper we begin the study of semisimple weak Hopf algebras
over an algebraically closed field of characteristic $0$. These are
the most accessible weak Hopf algebras that frequently appear in
applications. In particular, the representation-theoretic importance
of semisimple weak Hopf algebras can be seen from the result of V.~Ostrik \cite{O},
who proved that every semisimple rigid monoidal category with finitely
many classes of simple objects (such categories are called {\em fusion categories}) 
is equivalent to the representation category of a semisimple weak Hopf algebra.
Weak Hopf algebra techniques  were used to prove results about fusion categories
in \cite{N2, ENO}.
In \cite{ENO} it was shown that regular semisimple weak Hopf algebras 
and fusion categories 
do not admit continuous deformations (this extends the result of D.~Stefan
for semisimple Hopf algebras \cite{S}). Here the regularity condition
ensures the absence of ``trivial'' deformations \cite[Remark 3.7]{N2}.
In particular, this implies  that a classification of semisimple weak Hopf 
algebras is possible. Such a classification will include the classification
of fusion categories and module categories over them and will be useful 
in operator algebras, quantum field theory, and representation theory of quantum groups. 
In order to proceed with this classification project it is first  necessary to
bring the theory of semisimple weak Hopf algebras and fusion categories
to roughly the same level as that of semisimple Hopf algebras. 
The present paper is a step in this direction.

Semisimple Hopf algebras is one of the most
established subjects of Hopf algebra theory. A number of structural
and classification theorems is known, see surveys \cite{Mo1, Mo2, A}.
In this paper we extend several classical 
results from ordinary to weak Hopf algebras. 
Below we describe these results and organization of the paper.

In Section 2 we recall the definition and basic properties of weak Hopf
algebras and their representation categories. 

In Section 3 we first define
a canonical integral of a weak Hopf algebra $A$ and give a list of conditions
equivalent to the semisimplicity of $A$ in Proposition~\ref{a list}. Next,
for a semisimple $A$ we study two notions of dimensions of irreducible 
$A$-modules : the quantum and Frobenius-Perron dimensions
introduced in \cite{ENO}. We show that the vector of Frobenius-Perron dimensions
of irreducible $A$-modules is a Frobenius-Perron eigenvector (i.e., an eigenvector 
with strictly  positive entries) of the integer  matrix $\Lambda^t\Lambda$, where $\Lambda$ 
is the matrix describing the inclusion of the base algebra $A_s$ in $A$
(Proposition~\ref{LL}). This extends a result of \cite{BSz2} for weak Hopf $C^*$-algebras.
In Proposition~\ref{Ld} we give  in terms of $\Lambda$ a sufficient condition
of the triviality of a group-like element $G$ implementing the square of the 
antipode of $A$ (a pivotal element of $A$).  In the case when $A$ is pseudo-unitary
we show in Corollary~\ref{G=F} that $G=wS(w)^{-1}$ where $w$ is a Frobenius-Perron
element of the base $A_s$. In particular, 
in this case the eigenvalues of $S^2$ are all positive
and can easily be found.  

In Section 4 we study the Grothendieck ring of a semisimple weak Hopf algebra
and extend in Theorem~\ref{class eqn thm}
the Class Equation of G.~Kac and Y.~Zhu \cite{K, Z1}
(a version of this equation for pivotal fusion categories was proved in \cite{ENO}).

In Section 5 we establish analogues of the trace formulas of Larson and Radford 
\cite{LR1, LR2}.
The first formula (Corollary~\ref{TR S2}) expresses $\Tr(S^2|_A)$ in terms of dual integrals 
in $A$ and $A^*$. 
The second formula shows that the categorical dimension $\dim(A)$ of $A$ divides its
Frobenius-Perron dimension $\FPdim(A)$ in the ring of algebraic integers 
(Theorem~\ref{the 2nd trace}). As a consequence we obtain that  a weak Hopf algebra
is pseudo-unitary if and only if all the eigenvalues of $S^2$ are positive.

Section 6 deals with properties of module algebras and coalgebras over a semisimple
weak Hopf algebra $A$. In Theorem~\ref{J(M) is stable} we 
extend the result of \cite{Li} to prove that if $A$ is  a pseudo-unitary
weak Hopf algebra then the Jacobson radical
of any finite-dimensional
$A$-module algebra is $A$-stable. In Theorem~\ref{orbits} we extend the result
of  \cite{Z2} to Frobenius-Perron dimensions of comodule algebras over a semisimple
weak Hopf algebra $A$ : we show that if $M$
is an indecomposable finite-dimensional semisimple $A$-comodule algebra, then 
for all irreducible $M$-modules $M_1$ and $M_2$ the number
$\FPdim(A)\FPdim(M_1)\FPdim(M_2)/\FPdim(M)$, where 
$\FPdim$ denotes the Frobenius-Perron dimension,
is an algebraic integer. This generalizes the well-known 
fact that for a finite group $G$ the cardinality of a transitive $G$-set
divides the order of $G$.

\textbf{Acknowledgements.} This research  was supported by the NSF grant DMS-0200202.
The author is grateful to Pavel Etingof for poining out an error in an early version of
this paper and to Victor Ostrik and Leonid Vainerman
for valuable discussions. Thanks are also due to the referee whose comments helped
to improve the presentation.

\end{section}


\begin{section}
{Preliminaries}

\subsection{Definition and basic properties of weak Hopf algebras}

Let $k$ be an algebraically closed field of characteristic $0$.
All vector spaces and tensor products in this paper are over $k$.

One obtains a weak Hopf algebra by relaxing
the axioms related to the unit and counit
in the definition of an ordinary Hopf algebra.

\begin{definition}[\cite{BSz1, BNSz}] 
\label{finite weak Hopf algebra}
A {\em weak Hopf algebra} is a vector space $A$
with the structures of an associative algebra $(A,\,m,\,1)$ 
with a multiplication $m:A\otimes A\to A$ and unit $1\in A$ and a 
coassociative coalgebra $(A,\,\Delta,\,\epsilon)$ with a comultiplication
$\Delta:A\to A\otimes A$ and counit $\epsilon:A\to k$ such that:
\begin{enumerate}
\item[(i)] The comultiplication $\Delta$ 
is a (not necessarily unit-preserving) homomorphism of algebras:
\begin{equation}
\label{Delta m}
\Delta(gh) = \Delta(g) \Delta(h), 
\end{equation}
\item[(ii)] The  unit and counit satisfy the following identities: 
\begin{eqnarray}
\label{Delta 1}
(\Delta \otimes \id) \Delta(1) 
& =& (\Delta(1)\otimes 1)(1\otimes \Delta(1)) 
= (1\otimes \Delta(1))(\Delta(1)\otimes 1), \\
\label{eps m}
\epsilon(fgh) &=& \epsilon(fg\1)\, \epsilon(g\2h) = \epsilon(fg\2)\, 
\epsilon(g\1h),
\end{eqnarray}
\item[(iii)]
There is a linear map $S: A \to A$, called an {\em antipode}, such that
\begin{eqnarray}
m(\id \otimes S)\Delta(h) &=&(\epsilon\otimes\id)(\Delta(1)(h\otimes 1)),
\label{S epst} \\
m(S\otimes \id)\Delta(h) &=& (\id \otimes \epsilon)((1\otimes h)\Delta(1)),
\label{S epss} \\
S(h) &=& S(h\1)h\2 S(h\3),
\label{S id S}
\end{eqnarray}
for all $f,g,h\in A$.
\end{enumerate} 
\end{definition}

\begin{remark} We use Sweedler's notation for a comultiplication
in a coalgebra $C$, writing
$\Delta(c) = c\1 \otimes c\2$ for all $c\in C$.
\end{remark} 

Axioms (\ref{Delta 1})  and (\ref{eps m}) above 
weaken the usual bialgebra axioms 
requiring $\Delta$ to preserve the unit and
$\epsilon$ to be an algebra homomorphism.
Axioms (\ref{S epst}) and (\ref{S epss}) 
generalize the properties of the antipode in a Hopf algebra. 
The antipode $S$ of a weak Hopf algebra 
is an algebra and coalgebra antihomomorphism.
If $A$ is finite-dimensional then $S$ 
is bijective \cite[2.10]{BNSz}.

\begin{remark}
A weak Hopf algebra is a Hopf algebra if and only if 
$\Delta(1)=1\otimes 1$  and if and only if $\epsilon$ 
is a homomorphism of algebras.
\end{remark}

A {\em morphism} between weak Hopf algebras $A_1$ and $A_2$
is a map $\phi : A_1 \to A_2$ which is both algebra and coalgebra
homomorphism preserving $1$ and $\epsilon$
and which intertwines the antipodes of $A_1$ and $A_2$,
i.e., $ \phi\circ S_1 = S_2\circ \phi$. The image of a morphism is
clearly a weak Hopf algebra. 

For a finite-dimensional $A$
there is a natural weak Hopf algebra
structure on the dual vector space $A^*=\Hom_k(A,k)$ given by
\begin{eqnarray}
& & \la \phi\psi,\,h \ra = \la (\phi\otimes\psi),\, \Delta(h)\ra, \\
& & \la \Delta(\phi),\,g\otimes h\ra =  \la \phi,\,gh\ra, \\
& & \la S(\phi),\,h\ra = \la \phi,\,S(h)\ra,
\end{eqnarray}
for all $\phi,\psi \in A^*,\, g,h\in A$. The unit
of $A^*$ is $\epsilon$  and the counit is $\phi \mapsto \la \phi,\,1\ra$.

The linear maps defined by (\ref{S epst}) and (\ref{S epss})
are called {\em target\,} and {\em source counital maps}
and denoted $\eps_t$ and $\eps_s$ respectively :
\begin{equation}
\label{def counital maps}
\eps_t(h) = \epsilon(1\1 h)1\2, \qquad
\eps_s(h) = 1\1 \epsilon(h1\2),
\end{equation}
for all $h\in A$.

An algebra $M$ is a left {\em $A$-comodule algebra} if $M$ is a left
$A$-comodule via $\delta : M \to A \otimes M : m \mapsto m\I \otimes m\II$ 
such that 
\begin{equation}
\delta (mn) =  \delta(m) \delta(n), \qquad
\delta(1) = (\eps_s \otimes \id)\delta(1),
\end{equation}
for all $m,n\in M$.

An algebra $M$ is a left {\em $A$-module algebra} if $M$ is a left
$A$-module via $h\otimes m \mapsto h\cdot m$ such that 
\begin{equation}
h \cdot (mn) = (h\1\cdot m)(h\2 \cdot n),\qquad
h\cdot 1 = \eps_t(h),
\end{equation}
for all $h\in A$ and $m,n\in M$.

The definitions of right module and comodule algebras are similar.
\begin{remark}
\label{mod algebra}
If $A$ is finite-dimensional then a left $A$-comodule algebra $M$
is a right $A^*$-module algebra via
\begin{equation}
m\cdot \phi = \la \phi,\, m\I\ra  m\II, \qquad \phi\in A^*, m\in M.
\end{equation}
\end{remark}

\subsection{Basic properties of weak Hopf algebras}

The images of the counital maps defined in \eqref{def counital maps},
\begin{equation}
A_t = \eps_t(A),  \qquad A_s = \eps_s(A)
\end{equation}
are semisimple subalgebras of $A$, called target and source {\em bases} 
or {\em counital subalgebras} of $A$.  
These subalgebras commute with each other; moreover
\begin{eqnarray*}
A_t &=& \{(\phi\otimes \id)\Delta(1) \mid \phi\in A^* \} 
= \{ h\in A \mid \Delta(h) = \Delta(1)(h\otimes 1) \}, \\
A_s &=& \{(\id \otimes \phi)\Delta(1) \mid \phi\in A^* \}
= \{ h\in A \mid \Delta(h) = (1\otimes h)\Delta(1) \}.
\end{eqnarray*}

For any algebra $B$ we denote 
by $Z(B)$ the center of $B$.

If $p\neq 0$ is an idempotent in $A_t\cap A_s\cap Z(A)$, then
$A$ is the direct sum of weak Hopf algebras 
$pA$ and $(1-p)A$. Consequently, we say that $A$ is {\em indecomposable}
if $A_t\cap A_s\cap Z(A) = k1$.

Every weak Hopf algebra $A$ contains a canonical {\em minimal\,}
weak Hopf subalgebra $A_{{\rm min}}$ generated, as an algebra, 
by $A_s$ and $A_t$ \cite[Section 3]{N2}. 
Obviously, $A$ is an ordinary Hopf algebra
if and only if  $A_{{\rm min}} = k1$. Minimal weak Hopf algebras over $k$,
i.e., those for which  $A= A_{{\rm min}}$,
were completely classified in \cite[Proposition 3.4]{N2}.

The restriction of $S^2$ to $A_{{\rm min}}$ is always an inner automorphism
of $A_{{\rm min}}$, see \cite{N2}. 

\begin{note}
In what follows
we will consider only weak Hopf algebras satisfying 
the following natural property :
\begin{equation}
\label{regularity property}
S^2|_{A_{{\rm min}}} = \id.
\end{equation}
\end{note}

\begin{definition}
\label{def regular}
We will call a weak Hopf algebra satisfying \eqref{regularity property}
{\em regular}. 
\end{definition}

\begin{remark} It was shown in \cite[6.1]{NV2} that every weak Hopf algebra
can be obtained by twisting a  regular weak Hopf algebra
with the same algebra structure. 
\end{remark}

\subsection{Rigid monoidal categories and fusion categories}
\label{fusion categories}
Recall that a {\em monoidal category} consists of a category
$\mathcal{C}$, a tensor product bifunctor $\boxtimes : \C\times \C\to \C$,
a unit object $E$, and natural equivalences
\begin{equation}
\label{nat equiv}
\boxtimes(\id \times \boxtimes) \cong \boxtimes (\boxtimes\times\id),\quad
?\boxtimes E \cong \id,\quad E\boxtimes ? \cong \id,
\end{equation}
satisfying the pentagon and triangle diagrams \cite{Ma, BK}.

If $V$ is an object in $\C$ then a right dual of $V$ is an object
$V^*$ together with two morphisms, $b_V : E\to V\boxtimes V^*$ and
$d_V : V^* \boxtimes V \to E$, called coevaluation and evaluation,
such that
\begin{eqnarray}
\label{rigidity 1}
(\id_V\boxtimes d_V)(b_V \boxtimes \id_V) &=& \id_V, \\
\label{rigidity 2}
(d_{V} \boxtimes \id_{V^*})(\id_{V^*}\boxtimes b_V) &=& \id_{V^*}.
\end{eqnarray}
The definition of a left dual object is similar, see  \cite{BK}.
A monoidal category $\C$ is  called {\em rigid} if every object in $\C$
has right and left dual objects.

\begin{definition}
A {\em fusion category} is an Abelian semisimple rigid monoidal category $\C$
that has finitely many simple objects and is such that the unit object $E$ 
is simple and all $\mbox{Hom}$-spaces of $\C$ are finite-dimensional.
\end{definition}

A left {\em module category} over a monoidal category $\C$ is an Abelian category
$\mathcal{M}$ together with a bifunctor $\otimes : \C \times \M \to \M$ and
natural equivalences 
\begin{equation}
\label{nat equiv1}
\otimes(\id \times \otimes) \cong \otimes (\boxtimes\times\id),\quad
E\otimes ? \cong \id,
\end{equation}
satisfying pentagon and triangle axioms; see \cite{O} for details.

\subsection{Representation category of a weak Hopf algebra and reconstruction
of fusion categories}
\label{Rep category}

The category $\Rep(A)$ of finite-dimensional left $A$-modules 
has a natural structure of a rigid tensor category that we describe next, 
following \cite{NTV}. The tensor product of two $A$-modules $V$ and $W$ is
given by 
\begin{equation}
V \boxtimes W : = V \otimes_{A_t} W
=\{ x\in V\otimes_k W|\Delta(1)x =x\},
\end{equation}
where the right action of $A_t$ on $V$ is by 
$vz :=S(z)v,\, z\in A_t,\, v\in V$ and the $A$-module structure defined 
via $\Delta$. The tensor product of morphisms is defined in an obvious way.
The unit object $E$ of $\Rep(A)$ is the target counital 
subalgebra $A_t$ with the action 
$h\cdot z = \eps_t(hz)$ for all $h\in A,\, z\in A_t$.

The unit object of $\Rep(A)$ is irreducible if and only if bases of $A$ intersect 
trivially with the center of $A$, i.e., $Z(A)\cap A_t = k$.

\begin{definition}
\label{connected}
If $Z(A)\cap A_t = k$ we will say that $A$ is {\em connected}.
We will say that that $A$ is {\em coconnected}
if $A^*$ is connected and that $A$ is {\em biconnected} if it is both connected 
and coconnected.
\end{definition}

\begin{remark}
$A$ is coconnected if and only if  $A_s\cap A_t =k$ \cite[3.11]{N1}.
\end{remark}

If $V$ is an $A$-module then $V^*:=\Hom_k(V,\, k)$ is also an
$A$-module via 
\begin{equation}
\la h\cdot\phi,\, v\ra = \la \phi,\, S(h)\cdot v \ra,
\end{equation}
for all $h\in A,\, \phi\in A^*,\, v\in V$.

For any $V$ in $\Rep(A)$, we define evaluation and coevaluation
morphisms 
\begin{equation*}
d_V : V^*\boxtimes V \to A_t, \qquad b_V: A_t \to V \boxtimes V^*,
\end{equation*}
as follows. For $\sum_j\, \phi^j \otimes v_j\in V^* \boxtimes V$, set
\begin{equation}
\label{dV}
d_V(\sum_j\, \phi^j \otimes v_j)= \sum_j\, \phi^j(1\1\cdot v_j)1\2.
\end{equation}
Let $\{f_i\}$ and $\{\xi^i\}$ be bases of $V$ and $V^*$ 
dual to each other, then  $\sum_i\,f_i\otimes \xi^i$ does not
depend on choice of these bases. Set
\begin{equation}
\label{bV}
b_V(z) =  z\cdot (\sum_i\, f_i \otimes \xi^i).
\end{equation}
It was checked in  \cite{NTV} that $d_V$ and $b_V$ are well defined 
$A$-linear maps satisfying identities \eqref{rigidity 1} and
\eqref{rigidity 2}.
Thus, $\Rep(A)$ becomes a rigid monoidal category. Clearly,
$\Rep(A)$ is a fusion category if and only if $A$ is semisimple
and connected.

Let $\M$ be a semisimple left module category over a fusion category
$\C$ such that $\M$ has finitely many simple objects and all 
$\mbox{Hom}$-spaces of $\M$ are finite-dimensional. It was shown in \cite{O}
that there exists a biconnected semisimple weak Hopf algebra $A$ such that $\Rep(A)
\cong \C$ as fusion categories and $\Rep(A_t)\cong \M$ as module categories
over $\C$. Such a weak Hopf algebra  $A$ is not unique. In particular, the
base $A_t$ is defined up to a Morita equivalence, and one can always choose
$A_t$ to be a commutative algebra, in which case $A$ is regular in the sense
of Definition~\ref{def regular}. Thus, one can study fusion categories
and their module categories using weak Hopf algebra techniques.

\subsection{Grothendieck rings and modules and Frobenius-Perron dimensions}
\label{K0 and FP}

Let $A$ be a semisimple weak Hopf algebra. 
The {\em Grothendieck ring} $K_0(A)$ is defined as follows.
As an Abelian group it is generated by characters $\chi_V$
of finite-dimensional $A$-modules $V$ with the addition $\chi_U + \chi_V
=\chi_{U \oplus V}$ and multiplication $\chi_U\chi_V =\chi_{U\boxtimes V}$
for all $A$-modules $U$ and $V$. The unit is the character of the trivial
$A$-module. 
The map $\chi_V \mapsto \chi_{V^*}$ extends to an anti-multiplicative involution 
of $K_0(A)$ since $S^2$ is an inner automorphism of $A$.
Note that our definition does not differ from the usual one 
(see, e.g., \cite{Lo}) since in the semisimple case
we can identify finite-dimensional $A$-modules with their characters.

The ring $K_0(A)$ has a $\mathbb{Z}_+$-basis consisting of characters 
$\{\chi_1,\dots,\chi_n\}$ of irreducible  $A$-modules $\{V_1,\dots,V_n\}$.
For every $A$-module $V$ with the character $\chi_V$ the matrix $[\chi_V]$
of the left multiplication 
by $\chi_V$ in the basis $\{\chi_1,\dots,\chi_n\}$ has non-negative entries
and is not nilpotent. By the Frobenius-Perron theorem \cite{Ga} the matrix
$[\chi_V]$ has a real non-negative eigenvalue. The largest  among such
eigenvalues we will call
the {\em  Frobenius-Perron dimension} of $V$ and denote $\FPdim(V)$.
It was shown in \cite{ENO} that
the map $\phi : V \mapsto \FPdim(V)$ defines a ring homomorphism
$K_0(A) \to \mathbb{R}$. This  $\phi$  is a unique homomorphism 
$K_0(A) \to \mathbb{R}$ with the property $\phi(\chi_i)> 0$ for all 
$i=1,\dots n$.

\begin{definition}
Define the {\em character algebra} of $A$ to be 
$R(A)=K_0(A)\otimes_\mathbb{Z} \mathbb{C}$. 
If $k=\mathbb{C}$, then $R(A)$
may be regarded as a subalgebra of $A^*$.
\end{definition}

If $M$ is a finite-dimensional semisimple $A$-comodule algebra via
$\delta : M \to A\otimes M$ then the Grothendieck group $K_0(M)$
becomes a left $K_0(A)$-module (with non-negative integer 
structure constants) via
\begin{equation}
\la \chi\xi,\, m\ra = \la \chi\otimes \xi,\, \delta(m)\ra, 
\end{equation}
for all characters $\chi\in K_0(A),\, \xi\in K_0(M)$ and  $m\in M$.
Let $R(M)=K_0(M)\otimes_\mathbb{Z} \mathbb{C}$.

\begin{definition}
\label{indecomposable}
We will say that $M$ is {\em indecomposable} if it is not a direct sum
of two non-trivial $A$-comodule algebras. 
\end{definition}

Let $M$ be an indecomposable finite-dimensional semisimple $A$-comodule algebra.
Let $M_1,\dots,M_t$ be irreducible $M$-modules and let $\xi_1,\dots, \xi_t$
be their characters. It was shown in \cite{ENO} that there exist positive numbers
$\FPdim(M_1),\dots, \FPdim(M_t)$, defined up to a common positive scalar
multiple, and called {\em Frobenius-Perron dimensions}
of $M_1,\dots, M_t$, such that 
\begin{equation}
\xi_f = \sum_{i=1}^t\, \FPdim(M_i)\xi_i \in R(M)
\end{equation}
is an eigenvector  for every $A$-character $\chi_V$ with the eigenvalue $\FPdim(V)$:
\begin{equation}
\chi_V\xi_f = \FPdim(V)\xi_f.
\end{equation}
The assignment $\xi_i \mapsto \FPdim(M_i)$ extends to a $K_0(A)$-module homomorphism
$\psi: K_0(M)\to \mathbb{R}$  (where $\mathbb{R}$ is a $K_0(A)$-module via
$\chi_V a =\FPdim(V)a, \, a\in \mathbb{R}$) with $\psi(\xi_i)>0 $ for all $i=1,\dots t$.
The following statement will be used in the sequel.

\begin{lemma}
\label{unique homo}
The above map $\psi$ is  a unique, 
up to a positive scalar multiple, $K_0(A)$-module homomorphism
from $K_0(M)$ to $\mathbb{R}$ with the property $\psi(\xi_i)>0 $ for all $i=1,\dots t$.
\end{lemma}
\begin{proof}
Let $[\chi_V]$ be the matrix of multiplication by $\chi_V$ in the basis
$\{\xi_1,\dots, \xi_t\}$ and let $\vec{m}$ be the column vector 
with entries $\psi(\xi_i),\, i=1,\dots, t$. Then $[\chi_V]\vec{m} = 
\FPdim(V)\vec{m}$. Since $\vec{m}$ has strictly positive entries
it belongs to the Frobenius-Perron eigenspace  of  $[\chi_V]$ for
every $A$-module $V$. Choosing $V$ in such a way that   $[\chi_V]$
has strictly positive entries we conclude that $\vec{m}$ is
defined up to a scalar multiple (since in this case the Frobenius-Perron 
eigenspace  of  $[\chi_V]$ is $1$-dimensional).
\end{proof}

\end{section}

\begin{section}
{Integrals, semisimplicity, and dimension theory}

\subsection{Integrals in weak Hopf algebras}

The following notion of an integral in a weak Hopf algebra is a generalization
of that of an integral in a usual Hopf algebra.

\begin{definition}[\cite{BNSz}]
\label{integral}
A left (respectively, right) {\em integral} in 
a weak Hopf algebra $A$ is an element $\ell\in A$ 
(respectively, $r\in A$) such that 
\begin{equation}
h\ell =\eps_t(h)\ell, \qquad (\mbox{respectively, }rh = r\eps_s(h)) 
\qquad \mbox{ for all } h\in A. 
\end{equation}
\end{definition}
The space of left (respectively, right)
integrals in $A$ is a right (respectively, left) 
ideal of  $A$ of dimension 
$\dim_k(A_t)$. We will denote the spaces of left and right integrals
of $A$ by $\int_A^l$ and $\int_A^r$.

Any left integral $\lambda$ in $A^*$ satisfies the following
invariance property :
\begin{equation}
\label{left invariance}
g\1 \la \lambda,\,hg\2\ra  = S(h\1) \la \lambda,\,h\2g \ra,
\qquad g,h\in A.
\end{equation}

In what follows we use the Sweedler arrows $\actl$ and $\actr$
for the dual actions :
\begin{equation}
\label{Sweedler arrows}
\la h\actl\phi,\, g\ra = \la \phi,\, gh\ra
\quad \text { and } \quad
\la \phi\actr h,\, g\ra= \la \phi,\, hg\ra.
\end{equation}   
for all $g,\,h\in A,\phi\in A^*$.

Recall that a functional $\phi\in A^*$ is {\em non-degenerate}
if its composition with the multiplication
defines a non-degenerate bilinear form on $A$. 
Equivalently, $\phi$ is non-degenerate if
the linear map $h\mapsto (h \actl \phi)$
is injective. An integral (left or right) in a weak Hopf algebra
$A$ is called {\em non-degenerate} if
it defines a non-degenerate functional on $A^*$. A left integral $\ell$
is called {\em normalized} if $\eps_t(\ell)=1$. 

It was shown by P.~Vecsernyes \cite{V} that a finite-dimensional
weak Hopf algebra always 
possesses a non-degenerate left integral. 
In particular, a finite-dimensional weak Hopf algebra is a
Frobenius algebra (this extends the well-known result of Larson
and Sweedler for usual Hopf algebras). It also follows that
$\int_A^l$ is a free right $A_t$-module of rank $1$ and 
a free right $A_s$-module of rank $1$ for which any non-degenerate
left integral can be taken as a basis.

Maschke's theorem for weak Hopf algebras \cite[3.13]{BNSz}
states that a weak Hopf algebra $A$ is semisimple if and only if  $A$ is
separable, and if and only if  there exists a normalized left integral in $A$.
In particular, every semisimple weak Hopf algebra is finite-dimensional. 

For a finite-dimensional $A$ there is a useful notion of duality between 
non-degenerate
left integrals in $A$ and $A^*$ \cite[3.18]{BNSz}. If $\ell$
is a non-degenerate left integral in $A$ then there exists 
a unique $\lambda\in A^*$ such that
$\lambda\actl \ell = 1$. This $\lambda$ 
is a non-degenerate left integral in $A^*$. Moreover,
$\ell\actl \lambda =\epsilon$. Such a pair of non-degenerate integrals 
$(\ell,\,\lambda)$ is called a pair of {\em dual} integrals.

A weak Hopf algebra $A$ is called {\em unimodular} if it has a
$2$-sided non-degenerate integral. A semisimple weak Hopf algebra
is unimodular \cite{BNSz}.

Following Drinfeld, \cite{D}, we
define the algebra of {\em generalized characters} of $A$ to be
\begin{equation}
\label{gen characters}
O(A) = \{ \phi\in A^* \mid \la \phi,\, gh \ra = \la \phi, h S^2(g)\ra,
\quad g,h\in A \}.
\end{equation}

It was shown in \cite{BNSz} that a left integral $\ell\in \int_A^l$
is $S$-invariant if and only if its dual left integral
$\lambda\in \int_{A^*}^l$ is 
a generalized character.

\begin{proposition}
If $A$ is unimodular then
$\{ \ell \in \int_A^l \mid  \ell = S(\ell) \} \cong A_t \cap Z(A)$.
In particular if $A$ is connected unimodular then there exist 
unique up to a scalar multiple
non-degenerate $S$-invariant $\ell\in \int_A^l$  and non-degenerate
left integral $\lambda\in O(A)\cap \int_{A^*}^l$ dual to each other. 
\end{proposition}
\begin{proof}
Let $\ell\in A$ be  a two-sided non-degenerate integral, then
any other element of $\int_A^l$ is of the form $\ell' =\ell z,\,
z\in A_t$. It is easy to see that $S(\ell') =\ell'$ if and only if
$z \in A_t \cap Z(A)$.
\end{proof}

\begin{definition}
\label{canonical left integral}
A {\em canonical left integral} in $A^*$ is a functional $\lambda$
defined by 
\begin{equation}
\label{canonical intl}
\la\lambda,\, h\ra = \Tr(L_h\circ S^2|_A ),
\end{equation}
where $h\mapsto L_h$
is the left regular representation of $A$.
\end{definition}

That the above $\lambda$ is an integral was shown, e.g.,  in \cite{BNSz} and \cite{ENO}.
It follows that if $A$ is connected semisimple, then $\lambda$ is a basis of
$O(A)\cap \int_{A^*}^l$. 

The following definition was given in \cite{BNSz}.

\begin{definition}
\label{Haar intl}
A {\em Haar integral} in a weak Hopf algebra $A$ is a 
normalized $2$-sided integral $\ell$ in $A$. Such an integral is necessarily
unique and $S$-invariant.
\end{definition} 

Below we give a list of equivalent conditions characterizing
the semisimplicity of a biconnected weak Hopf algebra in terms of the
canonical and  Haar integrals.

\begin{proposition}
\label{a list}
Let $A$ be a biconnected weak Hopf
algebra and let $\lambda\in A^*$ be the canonical
integral defined by \eqref{canonical intl}.
Then the  following conditions are equivalent.
\begin{itemize}
\item[(i)] $\Tr(S^2|_A)\neq 0$,
\item[(ii)] $A$  is semisimple,
\item[(iii)] $A^*$  is semisimple,
\item[(iv)] There is a Haar integral in $A^*$,
\item[(v)] $\lambda$ is a non-degenerate functional on $A$.
\end{itemize}
\end{proposition}
\begin{proof}
(i) $\Rightarrow$ (ii) was proved in \cite[Corollary 6.6]{N2} and
(ii) $\Rightarrow$ (i) follows from \cite[Theorem 3.2]{ENO}.
Since $\Tr(S^2|_A) = \Tr(S^2|_{A^*})$ we also get 
(i) $\Leftrightarrow$ (iii).

To prove (iii) $\Leftrightarrow$ (iv) observe that
by \cite[Theorem 3.27]{BNSz} the existence of the Haar integral
in $A$ is equivalent to $A^*$ being semisimple with $S^2(\phi)=
\gamma\phi\gamma^{-1}$ for $\phi\in A^*$ and
\begin{equation}
\label{dim not zero}
\Tr_V(\gamma^{-1}) \neq 0
\end{equation}
for every irreducible $A^*$-module $V$.  But this condition
is satisfied in every  semisimple weak Hopf algebra, 
since the following sequence of $A^*$-module homomorphisms
\begin{equation*}
 A^*_t \stackrel{b_{V^*}}{\longrightarrow} V^* \boxtimes V^{**}  
\stackrel{\id_{V^*} \boxtimes \gamma^{-1}}{\longrightarrow}
V^* \boxtimes V \stackrel{d_{V}}{\longrightarrow} A^*_t
\end{equation*}
is precisely the multiplication by $\frac{\Tr_V(\gamma^{-1})}{\dim_k(A_t)}$
and is non-zero by semisimplicity of $A^*$.

The equivalence (iv) $\Leftrightarrow$ (v) was proved in 
\cite[Proposition 3.26(i)]{BNSz}.
\end{proof}

\begin{remark}
When $A$ is not biconnected there are refined criteria for semisimplicity,
cf.\ \cite[6.4]{N2} and \cite[4.10]{ENO}.
\end{remark}

Let $\{p_\alpha\}$ be the set of primitive idempotents in $Z(A_s)$.

\begin{proposition}
\label{counitals on lambda}
For any biconnected weak Hopf algebra we have
\begin{eqnarray}
\label{epstlambda}
\eps_t(\lambda) &=& \frac{\Tr(S^2|_A)}{\dim_k(A_s)}\, \epsilon, \\
\label{epsslambda}
\eps_s(\lambda) &=& \sum_\alpha \, \frac{\Tr(S^2|_{p_\alpha A})}{\dim_k(p_\alpha A_s)} 
\,(p_\alpha\actl \epsilon).
\end{eqnarray}
\end{proposition}
\begin{proof}
Let $\lambda$ be the canonical integral defined by \eqref{canonical intl}.
The restriction of $\lambda$ to $A_{\text{min}} = A_s A_t$ is an integral
in $(A_{\text{min}})^*$ and hence there exists $y\in A_s$ such that
\begin{equation}
\la \lambda,\, x \ra = \Tr(xy|_{A_{\text{min}}}),\quad x\in  A_{\text{min}}.
\end{equation}

In particular, for $x\in A_t$ we get
$\la \lambda,\, x \ra = \epsilon(x)\epsilon(y)$.
On the other hand, 
\begin{equation*}
\Tr(S^2|_A) = \la \lambda,\, 1\ra = \Tr(y|_{A_{\text{min}}}) =\epsilon(y) \dim_k(A_s),
\end{equation*}
whence \eqref{epstlambda} follows.

For $x\in A_s$ we have $\la \lambda,\, x \ra = \epsilon(xy) \dim_k(A_s)$,
so, $\epsilon_s(\lambda) = \dim_k(A_t) (y\actl \epsilon)$.
Let $y = \sum_\alpha\, y_\alpha p_\alpha$ for some scalars $y_\alpha$.
To find these scalars, note that $\la \lambda,\, p_\alpha \ra
= \Tr(S^2|_{p_\alpha A})$ and also that 
\begin{equation*}
\la \lambda,\, p_\alpha\ra =y_\alpha  \Tr(p_\alpha|_{A_{\text{min}}}) =
y_\alpha \dim_k(A_s) \dim_k(p_\alpha A_s).
\end{equation*}
Thus,
\begin{equation}
y = \sum_\alpha\, \frac{\Tr(S^2|_{p_\alpha A})}{ 
\dim_k(A_s) \dim_k(p_\alpha A_s)}\, p_\alpha,
\end{equation}
which implies \eqref{epsslambda}.
\end{proof}

\subsection{Group-like and pivotal elements in a weak Hopf algebra}

In \cite{N2} a  {\em group-like} element of $A$  was defined as an
invertible element $g\in A$ such that
$\Delta(g) = (g\otimes g)\Delta(1) = \Delta(1)(g\otimes g)$. 
Group-like elements of $A$ form a group $G(A)$
under multiplication. This group has a normal subgroup 
\begin{equation}
G_0(A) := G(A_{{\rm min}}) =\{ yS(y)^{-1}\mid y\in A_s\}
\end{equation} 
of {\em trivial} group-like elements. If $A$ is finite-dimensional,
the quotient group  $\widetilde{G}(A) = G(A)/ G_0(A)$ is finite.
It was shown in \cite{N2} that if 
$\ell\in A$ and $\lambda\in A^*$ is a dual pair of left integrals,
then there exist group-like elements $\alpha\in G(A^*)$ and $a\in G(A)$,
called {\em distinguished} group-like elements,
whose classes in $\widetilde{G}(A^*)$ and $\widetilde{G}(A)$
do not depend on the choice of $\ell$ and $\lambda$, such that
\begin{equation}
S(\ell) = \alpha\actl \ell \qquad \mbox{ and  } \qquad S(\lambda)=a \actl \lambda.
\end{equation} 
(Note that $\alpha$ and $a$ themselves depend on the choice of $\ell$ and $\lambda$).

The following result is  an analogue of
Radford's formula \cite{R} for usual Hopf algebras.

\begin{theorem}\label{S4} \cite[Theorem 5.13]{N2} One has
\begin{equation}
S^4(h) = a^{-1}(\alpha\actl h \actr \alpha^{-1})a. 
\end{equation}
for all $h\in A$.
\end{theorem}

\begin{remark}
$A$ is unimodular if the coset of $\alpha$ in $G(A^*)$ is trivial. 
\end{remark}

Recall that a {\em pivotal structure} on  a rigid monoidal category $\C$
is an isomorphism between monoidal functors $i: \mbox{Id} \to **$.
A weak Hopf algebra $A$ is {\em pivotal} if there is a group-like 
element $G\in A$ such that $S^2(h) = GhG^{-1}$ for all $h\in A$.
In other words, $A$ is pivotal if and only if $\Rep(A)$ is a
pivotal category.
The element $G$ is called a {\em pivotal element} of $A$.

\subsection{Dimension theory for weak Hopf algebras}

Let $A$ be  a semisimple weak Hopf algebra
and let $V$ be a finite-dimensional  $A$-module. 
Below we define quantum and Frobenius-Perron  dimensions
of $A$-modules. These dimensions coincide with
the $k$-vector space dimensions if and only if $A$ is a 
usual Hopf algebra.

Let $\C$ be a pivotal fusion category with
an isomorphism between monoidal functors $i: \mbox{Id} \to **$.
One can define the {\em quantum
dimension} of any object $V$ in $\C$ by
\begin{equation}
\label{quantum dim}
\dim(V) := \Tr_V(i) = d_{V^*} \circ (i \otimes \id_{V^*})\circ b_V.
\end{equation}
Here the right hand side belongs to $\End_\C(E) = k$.
Let  $\C =\Rep(A)$ for a pivotal weak Hopf algebra $A$ and let $V$ be
a finite-dimensional  $A$-module. If $G$ is a pivotal element of $A$ then
\eqref{quantum dim} becomes
\begin{equation}
\dim(V) = \frac{\Tr_V(G)}{\dim_k(A_t)}. 
\end{equation}
Note that we use $\dim$ for the quantum dimension and $\dim_k$
for $k$-vector space dimension. 

\begin{remark}
The {\em squared norm } $|V|^2$ of a simple object $V$ of $\C$
can be defined without assuming existence of a pivotal structure \cite{Mu, ENO}.
If $\C =\Rep(A)$ for a semisimple $A$, then $\End_k(V)$ can be identified
with a minimal $2$-sided ideal of $A$ and
\begin{equation}
\label{dim V sq}
|V|^2 = \frac{\Tr(S^2|_{\End_k(V)})}{\dim_k(A_t)^2}. 
\end{equation}
This number is an algebraic integer, which is non-negative if
$k=\mathbb{C}$ \cite[Theorem 2.3]{ENO}.
\end{remark}

\begin{remark}
For a pivotal $A$ one has $|V|^2= \dim(V)\dim(V^*)$ \cite{Mu, ENO}.
In general, $\dim(V)\neq \dim(V^*)$, and so  $|V|^2\neq \dim(V)^2$.
This is why we use the term ``squared norm'' rather than ``squared
dimension.''
\end{remark}


When $A$ is a semisimple weak Hopf algebra and $\{ V_j\}_{j=1}^n$ 
are all the irreducible $A$-modules,
the {\em categorical dimension} of
$\C =\Rep(A)$ is defined as
\begin{equation}
\dim(\C) = \sum_{j=1}^n\, |V_j|^2.
\end{equation}
We will call $\dim(A) = \dim(\Rep(A))$ the {\em categorical  dimension} of $A$.
It follows from \eqref{dim V sq} that for $k=\mathbb{C}$ one has
\begin{equation}
\label{dimA}
\dim(A) = \frac{\Tr(S^2|_A)}{\dim_k(A_t)^2}.
\end{equation}
It was shown in \cite{ENO}  that $\dim(A)$ is an algebraic integer 
$\geq 1$.

\begin{note}
From now on let $A$ be a semisimple biconnected weak Hopf algebra
and let $k=\mathbb{C}$ be the field of complex numbers (although
some of the results below remain valid over an arbitrary algebraically closed
field $k$ of characteristic $0$ we restrict our attention to $\mathbb{C}$
since it is convenient to regard dimensions of $A$-modules as elements of the
ground field).
\end{note}

Recall from Section~\ref{K0 and FP} that  the
{\em Frobenius-Perron dimension} $\FPdim(V)$ of a finite-dimensional $A$-module $V$
is defined as the largest positive eigenvalue of the matrix
of multiplication by $\chi_V$ in the basis 
$\{\chi_1,\dots,\chi_n\}$ of  $K_0(A)$, where
$\{\chi_1,\dots,\chi_n\}$ are characters of irreducible 
$A$-modules $\{V_1,\dots,V_n\}$. Clearly, $\FPdim(V)$ is a positive
algebraic integer.

The {\em Frobenius-Perron dimension} of $A$ is then defined \cite{ENO} as:
\begin{equation}
\FPdim(A) := \sum_{j=1}^n \, \FPdim(V_j)^2.
\end{equation}

\begin{remark}
It follows from \cite[Section 8.2]{ENO} that $\FPdim(A) =\FPdim(A^*)$.
If $B\subset A$ is a biconnected semisimple weak Hopf subalgebra 
of $A$ then $\frac{\FPdim(A)}{\FPdim(B)}$ is an algebraic integer.
We generalize this result in  Section~\ref{zhus}.
\end{remark}

Let us define
\begin{equation}
\label{rho}
\rho:=\sum_j\, \FPdim(V_j) \Tr_{V_j}\in R(A).
\end{equation}
This $\rho$ formally replaces the character of the regular representation.
We will call $\rho$ the {\em Frobenius-Perron character} of $A$.
If $V$ is an $A$-module with the character $\chi$
then $\chi\rho =\rho\chi =\FPdim(V)\rho$.  In particular,
$\mbox{FPdim}(A)^{-1}\rho$
is a minimal idempotent in $R(A)$.

A relation between categorical and Frobenius-Perron dimensions was found
in \cite[Section 8.3]{ENO}. For any simple module $V$ over a semisimple 
weak Hopf algebra $A$ one has $|V|^2 \leq \FPdim(V)^2$, and hence,
$\dim(A) \leq \FPdim(A)$. The ratio ${\dim(A)}/{\FPdim(A)}$
is an algebraic integer $\leq 1$. We explore this relation further
when we derive the second trace formula in Theorem~\ref{the 2nd trace}.

\begin{proposition}
\label{FP dimensions in At}
Let $A$ be a semisimple weak Hopf algebra.
There exists a unique, up to a scalar, element $w\in Z(A_s)$
that has strictly positive eigenvalues and  satisfies
\begin{equation}
\label{vector w}
\Tr_V(z w S(w)^{-1}) = \FPdim(V) \epsilon(z)
\end{equation}
for any $z\in A_t$ and any $A$-module $V$. In particular,
$w$ satisfies
\begin{equation}
\FPdim(V) = \frac{\Tr_V(w S(w)^{-1})}{\dim_k(A_s)}.
\end{equation}
\end{proposition}
\begin{proof}
The target counital subalgebra $A_t$ is a left $A$-comodule algebra 
via the comultiplication. Since $A^*$ is connected, $A_t$ is indecomposable
in the sense of Definition~\ref{indecomposable}.
As in Section~\ref{K0 and FP},
the Abelian group $K_0(A_t)$ is a left $K_0(A)$-module via
\begin{equation}
\la \chi\xi,\, z\ra = \la \chi\otimes \xi,\, \Delta(z)\ra, \qquad z\in A_t, 
\end{equation}
for any characters $\chi$ of $A$ and $\xi$ of $A_t$. This module
structure naturally extends to the algebra $R(A_t) = K_0(A_t)\otimes_\mathbb{Z} \mathbb{C}$. 
It follows that there exists a unique, up to a scalar, $\xi_f\in R(A_t)$ with strictly
positive coordinates 
in the basis of characters $\{\xi_1,\dots,\xi_t\}$ of 
irreducible $A_t$-modules such that
\begin{equation}
\label{hiksi}
\Tr_V\xi_f = \FPdim(V) \xi_f,  
\end{equation}
for any $A$-module $V$ with character $\Tr_V$. Let $w=\xi_f\actl 1\in Z(A_s)$, then
\begin{equation}
\Tr_V \actl w = \FPdim(V)w.
\end{equation}
Thus, $w$ is a Frobenius-Perron eigenvector of the linear transformation
$\Tr_V\actl$ of $Z(A_s)$ with the eigenvalue  $\FPdim(V)$. Note that
the corresponding eigenspace is $1$-dimensional by the Frobenius theorem \cite{Ga}.
The last equation, in turn, is equivalent to
\begin{equation}
\label{qf}
1\1 w^{-1} \Tr_V(1\2 w ) = \FPdim(V), 
\end{equation}
which is equivalent to \eqref{vector w} (this can be seen by applying
$z\actl \epsilon$ to both sides). Therefore,  $w$ is unique
up to a scalar multiple and has positive eigenvalues.
\end{proof}

\begin{definition}
\label{FPdim in As}
A non-zero  element $w\in Z(A_s)$ 
will be called a {\em Frobenius-Perron element} of $A_s$
if it has positive eigenvalues and  satisfies
\begin{equation} 
\label{f is PF}
\Tr_V \actl w = \FPdim(V)w
\end{equation}
for all finite-dimensional  $A$-modules $V$. This element 
is defined up to a non-zero scalar multiple.
\end{definition}

\begin{corollary}
\label{rho actl 1}
Every Frobenius-Perron element of $A_s$ is a multiple of $\rho\actl 1$.
\end{corollary}

Let $A$ be a semisimple weak Hopf algebra.
To fix the notation in what follows we let
\begin{equation}
A \cong \oplus_{j=1}^n \, \End_k(V_j), \qquad 
A_s\cong \oplus_{\alpha=1}^l\, M_{n_\alpha}(\mathbb{C}),
\end{equation}
where $M_n(\mathbb{C})$ is the algebra of $(m\times m)$-matrices  over $\mathbb{C}$.
Let $\Lambda =(\Lambda_{\alpha j})$ be the $(l\times n)$ matrix 
of the inclusion $A_s\subset A$, i.e.,  $\Lambda_{\alpha j}$ is the 
multiplicity of the simple $A_s$-module corresponding to $M_{n_\alpha}(\mathbb{C})$
in the restriction of the $A$-module $V_j$ to $A_s$.

\begin{definition}
\label{index}
The number 
\begin{equation}
\mu(A) =\sum_{j=1}^n\, \dim_k(V_j) \FPdim(V_j)
\end{equation}
will be called an {\em index} of $A$, cf.\ \cite[4.5]{BSz2}.
\end{definition}
In other words, $\mu(A)$ is the Frobenius-Perron dimension of 
the regular representation of $A$.
The name index comes from the subfactor theory \cite{GHJ}:  if a 
subfactor comes from the crossed product with a weak Hopf $C^*$-algebra $A$ 
then the Jones index of the subfactor is equal to the index of $A$.

Define the {\em Frobenius-Perron dimension vector} of $A$ as 
\begin{equation}
\vec{f} =(\FPdim(V_j))_{j=1}^n.
\end{equation}
If $A$ is pivotal with a pivotal element $G$ define the {\em quantum
dimension vector} of $A$ as
\begin{equation}
\vec{d} =(\dim(V_j))_{j=1}^n.
\end{equation}
Recall from \eqref{rho} the canonical character 
$\rho= \sum_j\, \FPdim(V_j)\Tr_{V_j}$ of $A$. 
Let $p_\alpha,\, \alpha=1,\dots, l$ be primitive central idempotents
of $A_s$. 

Let $w_\alpha = \frac{\la \rho,\,p_\alpha \ra}{n_\alpha^2 \mu(A)}$
and $v_\alpha =w_\alpha n_\alpha$.
It follows from Corollary~\ref{rho actl 1} that 
\begin{equation*}
 w = \sum_\alpha\,
w_\alpha p_\alpha
= \frac{1}{\mu(A)}\, \rho\actl 1
\end{equation*}
is a Frobenius-Perron element of $A_s$. Let
\begin{equation}
\vec{v} = (v_\alpha)_{\alpha=1}^l =  \left( \frac{\la \rho,\,p_\alpha \ra}
{n_\alpha \mu(A)}\right)_{\alpha=1}^l.
\end{equation} 
Note that if $S^2=\id$ then $v_\alpha=n_\alpha$, 
i.e., $\vec{v}$ is the dimension vector of $A_s$.

\begin{remark}
The inclusion matrix $\Lambda$, vector $\vec{v}$,  
and index $\mu(A)$ are not invariants
of the representation category of $A$, i.e., there exist semisimple
weak Hopf algebras with monoidally equivalent representation categories but
different $\Lambda$, $\vec{v}$, and $\mu(A)$.
\end{remark}

\begin{proposition}
\label{LL}
Let $A$ be a semisimple weak Hopf algebra.
Then $\Lambda \vec{f} = \mu(A)\vec{v}$
and $\Lambda^t \vec{v} = \vec{f}$, where $\Lambda^t$ is
the transpose of $\Lambda$. Consequently,
\begin{equation}
\Lambda^t\Lambda \vec{f} =\mu(A) \vec{f} 
\qquad \mbox{ and }\qquad
\Lambda\Lambda^t \vec{v} =\mu(A) \vec{v}.
\end{equation}
\end{proposition}
\begin{proof}
We have $\FPdim(V^*)=\FPdim(V)$ for any simple $A$-module $V$. 
Let $f_j,\, j=1,\dots, n$ and $v_\alpha,\, \alpha=1,\dots,l$
be the coordinates of vectors $\vec{f}$ and $\vec{v}$.
We compute
\begin{eqnarray*}
\sum_j\, \Lambda_{\alpha j} f_j 
&=& \sum_j\, \frac{\Tr_{V_j}(p_\alpha)}{n_\alpha} \FPdim(V_j)\\
&=& \frac{\la \rho,\,p_\alpha\ra}{n_\alpha} = \mu(A)v_\alpha,\\
\sum_\alpha\, \Lambda_{\alpha j}  v_\alpha
&=& \sum_\alpha\, \frac{\Tr_{V_j}(p_\alpha)}{n_\alpha} 
    \frac{\la \rho,\,p_\alpha \ra}{n_\alpha \mu(A)} \\
&=& \mu(A)^{-1} \Tr_{V_j^*}(1\2) \la \rho, 1\1 \ra \\
&=& \mu(A)^{-1} \la \rho \Tr_{V_j^*}, 1\ra \\
&=& \FPdim(V_j) \mu(A)^{-1} \la \rho,\, 1\ra = f_j,
\end{eqnarray*}
whence the statement of the proposition follows.
\end{proof}

\begin{remark}
A version of Proposition~\ref{LL} was proved in \cite[Theorem 4.5]{BSz2}
for weak Hopf $C^*$-algebras.
\end{remark}

\begin{remark}
It follows from Proposition~\ref{LL} that $\rho$ is a Markov
trace \cite{GHJ} for the inclusion $A_s\subset A$.
\end{remark}

The Frobenius-Perron vectors of $A_s$ and $A^*_s$ turn out to
be colinear as we show in the next Proposition. 
In the language of \cite{ENO}
this means that the vectors of Frobenius-Perron dimensions of
a finite semisimple module category $\mathcal{M}$ viewed as a module category
over a fusion category $\C$ and over its dual $\C^*_\M$
are proportional.

\begin{proposition}
\label{w and w*}
Let  ${w}_A$ and ${w}_{A^*}$ be the Frobenius-Perron
 elements of $A_s$ and $A^*_s$ respectively.  Then
${w}_{A^*}$ and $w_A\actl \epsilon $ are scalar
multiples of one another.
\end{proposition}
\begin{proof}
Recall that ${w}_A$ and ${w}_{A^*}$ are defined 
by \eqref{f is PF}
in terms of the left coregular actions of $R(A)$ and $R(A^*)$
on $Z(A_s)$ and $Z(A^*_s)$ respectively. 
If we identify
$Z(A^*_s)$ with $Z(A_s)$ via $\phi\mapsto (\phi\actl 1),\,
\phi\in Z(A_s^*)$, then the above actions make $Z(A_s)$
an $R(A)-R(A^*)$ bimodule with $w_A$ being an eigenvector
of all $\chi\in R(A)$ and $x\in R(A^*)$. Since $w_{A^*}$ can be 
normalized to have positive entries it follows
that it belongs to the Frobenius-Perron eigenspace 
of $R(A)$ that also contains  $w_{A}\actl \epsilon$.
\end{proof}

\subsection{Implementation of the antipode and pseudo-unitary weak Hopf algebras}

\begin{definition}
\label{pseudo-unitary WHA}
A semisimple weak Hopf algebra $A$ is said to be {\em pseudo-unitary}
if $\dim(A) = \FPdim(A)$.
\end{definition}

Any semisimple Hopf algebra $A$ is pseudo-unitary since in this
case  $\dim(A) = \FPdim(A) =\dim_\mathbb{C}(A)$.
Weak Hopf $C^*$-algebras considered in \cite{BNSz, BSz1, BSz2, N1} are also examples
of pseudo-unitary weak Hopf algebras. However, there are many
semisimple weak Hopf algebras that are not pseudo-unitary \cite[Section 8]{ENO}.

For a pseudo-unitary $A$ one has  $|V|^2 =\FPdim(V)^2$ for all
irreducible $A$-modules $V$. It was shown in \cite[Section 8.4]{ENO}
that a pseudo-unitary $A$ has a unique pivotal element $G$ with respect 
to which the quantum dimensions of simple objects coincide with their
Frobenius-Perron dimensions. We will call this $G$ a {\em canonical
pivotal element} of $A$. Below we prove that
the canonical pivotal element is a trivial group-like element with
positive eigenvalues. The next Proposition establishes a sufficient
condition of the triviality of a pivotal element.

\begin{proposition}
\label{Ld}
Let $A$ be a semisimple weak Hopf algebra that has a pivotal element $G$.
Let $\Lambda$ be the inclusion matrix of  $A_s\subset A$ 
and let $\vec{d}=(\dim(V_j))_{j=1}^n$ be the vector
of quantum dimensions of irreducible $A$-modules.
If $\Lambda \vec{d}\neq \vec{0}$ the
$G$ is a trivial group-like element.
\end{proposition}
\begin{proof}
The left coregular action
\begin{equation}
\phi\mapsto (G\actl\phi), \qquad \phi\in A^*,
\end{equation}
of a group-like element $G\in A$ on $A^*$ is an algebra automorphism of $A^*$.
This  automorphism preserves the minimal $2$-sided ideal $I=\End_k(A^*_t)$ 
of $A^*$ corresponding to the trivial $A^*$-module
if and only if $G$ is trivial \cite{N2}.

The canonical integral $\lambda\in A^*$ defined by \eqref{canonical intl}
is a rank $1$ element of $I$. Therefore, $G$ is a trivial 
group-like element if and only if $(G^{-1}\actl\lambda)\in I$, i.e., if
$G^{-1}\actl\lambda$ acts as non-zero in the trivial representation of $A^*$. 
Since $G^{-1}\actl\lambda = \sum_i\, \Tr_{V_i}(G^{-1})\Tr_{V_i}$, the last
condition is equivalent to
\begin{equation}
\sum_i\, \Tr_{V_i}(G)\Tr_{V_i}|_{A_tA_s} \neq 0.
\end{equation}
In particular, this condition follows if
$\sum_i\, \Tr_{V_i}(G)\Tr_{V_i}(p_\alpha)\neq 0$
for some primitive idempotent $p_\alpha\in Z(A_s)$. Since
$\Lambda_{\alpha i} = \frac{\Tr_{V_i}(p_\alpha)}{n_\alpha}$ and
$d_i = \frac{\Tr_{V_i}(G)}{\dim_k(A_s)}$ we get the result.
\end{proof}

\begin{corollary}
\label{positive triviality}
Suppose that $A$ has a pivotal element $G$ with respect to which
all quantum dimensions of simple $A$-modules are non-negative.
Then $G$ is a trivial group-like element.
\end{corollary}

\begin{corollary}
\label{pseudo-unital triviality}
Let $A$ be a pseudo-unitary weak Hopf algebra and let $G\in A$ be the pivotal
element with respect to which  the quantum dimensions of
all simple objects coincide with their Frobenius-Perron dimensions.
Then $G$ is a trivial group-like element.
\end{corollary}

\begin{corollary}
\label{G=F}
If $A$ is pseudo-unitary, then $G=w S(w)^{-1}$, where $w\in Z(A_s)$
is a Frobenius-Perron element of $A_s$, is the canonical
group-like element of $A$. 
\end{corollary}
\begin{proof}
Up to a  scalar multiple we have $w=\rho\actl 1$. Since $A$
is pseudo-unitary, its canonical pivotal
element $G$ is trivial by Corollary~\ref{pseudo-unital triviality},
i.e., $G= gS(g)^{-1}$ for some $g\in Z(A_s)$.
We also have $\rho= \dim_k(A_s)^{-1}(G^{-1}\actl \lambda)$, 
where $\lambda$ is the canonical integral in $A^*$. We compute
\begin{eqnarray*}
\dim_k(A_s) w
&=& 1\1 \la \lambda,\, G^{-1}1\2 \ra \\
&=& g1\1\la \lambda,\, g^{-1}1\2\ra \\
&=& g (\lambda \actl g^{-1}).
\end{eqnarray*}
Since $\lambda$ is an integral and $A$ is biconnected,
we have $\lambda \actl g^{-1} \in A_s\cap A_t = k1$ is
a scalar, hence $w$ and $g$ are proportional and $wS(w)^{-1} 
= gS(g)^{-1} = G$.
\end{proof}

\begin{remark}
It follows that if $A$ is pseudo-unitary then the eigenvalues of $S^2$ 
belong to the set $\left\{\frac{w_\alpha w_\beta}{w_\gamma w_\delta}
\mid \alpha,\beta,\gamma,\delta =1,\dots, l\right\}$,
where $w$ is a Frobenius-Perron element of $A_s$.  
In particular, the eigenvalues of $S^2$
of a pseudo-unitary weak Hopf algebra 
are strictly positive.
Below we will see  (Corollary~\ref{ps.u = positivity}) that the converse
is also true.
\end{remark}

\end{section}


\begin{section}
{The class equation}

\subsection{The character algebra and the Grothendieck ring of
a weak Hopf algebra}

Let $A$ be a semisimple weak Hopf algebra, $K_0(A)$ be the Grothendieck
ring of $A$ defined in Section~\ref{K0 and FP}, and $R(A) \subset A^*$
be the character algebra.

\begin{remark}
Note that, in general, $\epsilon\not\in R(A)$. 
The identity element of $R(A)$ is the character of the trivial
representation of $A$, $\chi_1 = \epsilon\1\epsilon\2$. 
\end{remark}

In what follows we will identify $K_0(A)$ with a subring of $R(A)$.
For any finite-dimensional $A$-module $V$ let $\chi_V$ be its character.
The operation of taking the dual module gives rise to an
anti-isomorphism of $R(A)$ :
\begin{equation}
\chi_V\mapsto \chi_{V^*} = \chi_V\circ S,
\end{equation}
for any character $\chi_V$. The map $*$ extends to an involutive anti-linear
algebra anti-homomorphism.

Let $V_1,\dots,V_n$ be a complete
set of simple $A$-modules, with $V_1 =A_t$ the trivial $A$-module.
Let $\chi_1,\dots,\chi_n$ be the corresponding characters
and let $\chi_j^* = \chi_j\circ S$ be the character of $V_j^*$.

The following lemma is standard.

\begin{lemma}
\label{bilinear}
The bilinear form $( \cdot\,,\, \cdot) :
R(A)\times R(A) \to \mathbb{C}$ defined by
\begin{equation}
( \chi_V,\, \chi_W) :=\dim_k \Hom_A(V,\,W^*)
\end{equation}
is non-degenerate, associative, $*$-invariant, and symmetric.
\end{lemma}
\begin{proof}
That the given form is non-degenerate follows from the fact
that the bases $\{\chi_1,\dots,\chi_n \}$ and $\{\chi_1^*,\dots,\chi_n^* \}$
of $R(A)$ are dual to each other with respect to  $( \cdot, \cdot)$, 
i.e., $( \chi_i,\, \chi_j^*) =\delta_{ij}$ for all $i$ and $j$. 
The associativity of  $( \cdot, \cdot)$ follows from the isomorphism
of vector spaces
\begin{equation*}
\Hom_A(V\boxtimes U,\, W^*) \cong \Hom_A (V\boxtimes U
\boxtimes W,\, V_1)
\cong \Hom_A(V,\, W^* \boxtimes U^*),
\end{equation*}
for all $A$-modules $V,\,U,\,W$.
Clearly, the form is $*$-invariant and symmetric.
\end{proof}

\begin{lemma}
\label{R(A) is ss}
The character algebra $R(A)$ is semisimple.
\end{lemma}
\begin{proof}
We claim that $\sum_i\,\chi\chi_i \otimes \chi_i^*
= \sum_i\,\chi_i \otimes \chi_i^*\chi$ for all $\chi\in R(A)$.
Indeed, evaluating both sides of this equality 
against $\chi_k\otimes\chi_l$ we get $(\chi,\,\chi_l\chi_k)$,
whence the claim follows by Lemma~\ref{bilinear}.

Next, $( \sum_i\,\chi_i\chi_i^*\chi,\,\chi^*) = 
\sum_i\, (\chi_i\chi,\,(\chi_i\chi)^* )\geq 0$,
where the equality occurs only for $\chi=0$. This implies that
$(\sum_i\,\chi_i\chi_i^*)\chi\neq 0$ for $\chi\neq 0$,
and so $\sum_i\,\chi_i\chi_i^*$ is an invertible
central element in $R(A)$.
This means that $R(A)$ is separable, and therefore, semisimple.
\end{proof}

\begin{remark}
The element $\sum_i\,\chi_i\chi_i^*\in K_0(A)$ is the character
of the adjoint representation of $A$ on the centralizer
$C_A(A_t)$ of its target base, cf.\ \cite{NTV}. 
\end{remark}

\subsection{The class equation for weak Hopf algebras} 

We extend the argument of \cite{Lo} to obtain an analogue
of the class equation of Kac \cite{K} and Zhu \cite{Z1} showing
that the categorical dimension $\dim(A)$ of a semisimple weak
Hopf algebra $A$ is equal to a sum of its divisors in the 
ring of algebraic integers of $\mathbb{C}$.

Let $e\in R(A)$ be a primitive  idempotent in $R(A)$, and $\hat{e}$ be its 
central support in $R(A)$, i.e., the primitive idempotent in $Z(R(A))$ 
such that $e\hat{e}=e$. Let $m$ be the dimension of the simple $R(A)$-module
$R(A)e$, then $R(A)\hat{e} \cong (R(A)e)^m$ as $R(A)$-modules.

Let $\mu$ be the character of $R(A)\hat{e}$ and 
let $\omega = \frac{1}{m}\mu$.

\begin{lemma}
\label{hat e} We have
$\hat{e} = m \mu(\tilde{e})^{-1} \tilde{e} =\omega(\tilde{e})^{-1}  \tilde{e}$, 
where $\tilde{e} =\sum_i\,\mu(\chi_i^*)\chi_i$ is an invertible element in $Z(R(A))$.
\end{lemma}
\begin{proof}
That $\tilde{e}$ is central and invertible  follows from the proof of Lemma~\ref{R(A) is ss}. 
Clearly, $\tilde{e}e'=0$ for all primitive idempotents
$e'\in Z(R(A))$ such that $e'\neq \hat{e}$, therefore 
$\tilde{e}$  is proportional to $\hat{e}$ and
$\hat{e} = m \chi(\tilde{e})^{-1} \tilde{e}$.
\end{proof}

Let $d= \dim_k(A_t)$.

\begin{lemma} 
\label{1112}
In any finite-dimensional weak Hopf algebra $A$ we have
\begin{equation}
\Tr(S^2|_{1\11\2A})=\frac{1}{d} \Tr(S^2|_A).
\end{equation}
\end{lemma}
\begin{proof}
Let $\lambda\in A^*$ be the canonical left integral defined in 
\eqref{canonical intl}. Let $y\in Z(A_s)$ be as in 
the proof of Proposition~\ref{counitals on lambda}, i.e., such that
\begin{equation}
\label{integral restricted}
\la \lambda,\, x\ra = \Tr(xy|_{A_{\text{min}}}),\qquad 
\forall x\in A_{\text{min}}.
\end{equation}
Observe that 
\begin{equation}
\Tr(zy|_{A_{\text{min}}}) =\epsilon(z)\epsilon(y),\qquad \forall z\in A_t,\, y\in A_s.
\end{equation}
Taking $x=1\11\2$ in \eqref{integral restricted} we get
\begin{eqnarray*}
\Tr(S^2|_{1\11\2A})
&=& \Tr(1\11\2y|_{A_{\text{min}}}) \\
&=& \epsilon(1\1y)\epsilon(1\2) =\epsilon(y) \\
&=& \frac{1}{d}  \Tr(y|_{A_{\text{min}}}) = \frac{1}{d} \Tr(S^2|_A).
\end{eqnarray*}
\end{proof}

\begin{proposition}
\label{we have}
We have $\Tr(S^2|_{eA^*})\neq 0$ and
\begin{equation}
 \dfrac{\Tr(S^2|_A)}{d\, \Tr(S^2|_{eA^*})} =\omega(\tilde{e}).
\end{equation}
\end{proposition}
\begin{proof}
Let $\ell$ be the left canonical integral in $A$, cf. \eqref{canonical intl}, i.e.,
\begin{equation}
\label{ell}
\la \phi,\, \ell \ra = \Tr(L_\phi\circ S^2|_{A^*}),\qquad \forall \phi\in A^*,
\end{equation}
where $\phi\to L_\phi$ is the left regular representation of $A^*$.
From Lemma~\ref{hat e} we get
\begin{equation*}
\Tr(S^2|_{\hat{e}A^*}) = \la \hat{e},\,\ell \ra 
= \omega(\tilde{e})^{-1} \la \tilde{e},\, \ell\ra. 
\end{equation*}
On the other hand, since $\la \chi_i,\, \ell\ra =0$ if $i\neq 1$, we have
\begin{eqnarray*}
\la \tilde{e},\, \ell \ra
&=& \sum_i \mu(\chi_i^*) \la \chi_i,\, \ell \ra \\
&=& \mu(\chi_1) \la \chi_1,\, \ell \ra \\
&=& m \Tr(S^2|_{\epsilon\1\epsilon\2A^*}).
\end{eqnarray*}
Combining two last formulas we get
\begin{equation*}
\omega(\tilde{e}) \Tr(S^2|_{\hat{e}A^*})=
m\Tr(S^2|_{\epsilon\1\epsilon\2 A^*}).
\end{equation*}
Next, since $S^2|_{R(A)} =\id$ the value $\Tr(S^2|_{{e}A^*})$
is the same for all primitive idempotents $e\in \hat{e}R(A)$,
hence  $\Tr(S^2|_{\hat{e}A^*}) =m \Tr(S^2|_{eA^*})$. Using this
and Lemma~\ref{1112} we obtain the result.
\end{proof}

Recall from equation~\eqref{dimA} that the categorical dimension of 
a semisimple weak Hopf algebra $A$ is
\begin{equation*}
\dim(A) =\frac{\Tr(S^2|_A)}{d^2} = \frac{\Tr(S^2|_{A^*})}{d^2}.
\end{equation*}

\begin{theorem}[The Class Equation]
\label{class eqn thm}
Let $e_1,\dots e_k$ be primitive idempotents of $R(A)$ such that 
$e_ie_j =\delta_{i,j}e_i$ and $\sum_i\,e_i =\epsilon\1\epsilon\2$. Then
\begin{equation}
\label{class eqn}
\dim(A) = \sum_i\, \frac{\Tr(S^2|_{e_iA^*})}{d},
\end{equation}
where the numbers $n_i=\dfrac{\Tr(S^2|_{e_iA^*})}{d},\, i=1,\dots, k$ 
are such that the ratios $\dim(A)/n_i$ are 
algebraic integers in $\mathbb{C}$.
\end{theorem}
\begin{proof}
In view of Proposition~\ref{we have} it only remains to show that
$\omega(\tilde{e})$ is an algebraic integer. We reproduce here the argument
of \cite{Lo} for the sake of completeness. Note first that for any
character $\psi$ of $R(A)$ we have $\psi(K_0(A)) \subset \algint(\mathbb{C})$. 
Indeed, the regular representation of $K_0(A)$ on itself is faithful, and so
each $\chi_i\in K_0(A)$ satisfies a monic polynomial over $\mathbb{Z}$, e.g.,
the characteristic polynomial of the matrix $(N_i)_{jk} = 
(\chi_i\chi_j,\,\chi_k^*)$. Hence, so does the image $\pi(\chi_i)$
under any representation $\pi$ of $R(H)$. Therefore, all eigenvalues
of $\pi(\chi_i)$ and $\Tr(\pi(\chi_i))$ belong to $\algint(\mathbb{C})$.

This means that $\tilde{e} = \sum_i\,\mu(\chi_i^*)\chi_i$ is a linear
combination of elements of $K_0(A)$ with coefficients from $\algint(\mathbb{C})$. 
Since $\omega : Z(R(A))\to \mathbb{C}$ is an algebra  homomorphism, $\omega(\chi_i)
\in \algint(\mathbb{C}),\, i=1,\dots, n$. Thus, $\omega(\tilde{e})\in \algint(\mathbb{C})$.
\end{proof}

\begin{remark}
It was shown in \cite{ENO} that if $A$ is pivotal then the numbers $n_i$
are algebraic integers equal to the quantum  dimensions of irreducible submodules
of the $D(A)$-module induced from the trivial $A$-module, where $D(A)$ is
the Drinfeld double of $A$ \cite{NTV}. 
\end{remark}

\end{section}


\begin{section}
{Trace formulas}
\label{chapter 5}

Recall that in the case of a usual semisimple Hopf algebra $A$ 
at least one of the summands in the Class Equation 
\eqref{class eqn} equals $1$. Namely, 
if $\Tr$ denotes the trace of the left regular representation
of $A$, then $\rho=\dim_\mathbb{C}(A)^{-1}\Tr$ is a minimal idempotent in $R(A)$
and $\Tr(S^2|_{\rho A^*})= 1$. Below we extend this observation
to weak Hopf algebras and also establish an analogue of the trace
formula \cite[Formula~(6)]{LR2}.

\subsection{The first trace formula}

A weak Hopf algebra analogue of the Larson-Radford formula
for $\Tr(S^2)$ \cite{LR2} was established in \cite{N2}. 
As in the case of usual Hopf algebras, 
this formula is related to the semisimplicity of 
the corresponding weak Hopf algebra and its dual.

It follows from \eqref{left invariance} that for a non-degenerate
integral $\ell$ in $A$ with the dual integral $\lambda$, the element
$(\ell\2 \actl \lambda) \otimes S^{-1}(\ell\1) \in A^*\otimes A$ 
is the dual bases tensor. This implies that for every $T\in \End_\mathbb{C}(A)$
we have
\begin{equation}
\label{Tr T}
\Tr(T) = \la \lambda, T(S^{-1}(\ell\1))\ell\2).
\end{equation}
In particular, for $T = S^2$, we get the following 
analogue of \cite[Theorem~2.5(a)]{LR2}
with counits replaced by counital maps.

\begin{corollary}[The first trace formula \cite{N2}]
\label{TR S2}
Let $\ell\in A$ be a non-degenerate integral in a finite-dimensional 
weak Hopf algebra $A$ and let $\lambda$ be the dual integral. Then
\begin{equation}
\label{Formula TR S2}
\Tr(S^2|_A) = \la \eps_s(\lambda), \eps_s(\ell)\ra.
\end{equation}
\end{corollary}

\subsection{The second trace formula}
We will derive a relation between the categorical and Frobenius-Perron dimensions
of a semisimple weak Hopf algebra $A$ that extends the result of \cite{LR2}.

Recall that $\chi_1=\epsilon\1\epsilon\2$ is the character of the trivial 
left $A$-module $A_t$.
\begin{lemma}
\label{chi0}
$\Tr(S^2|_{\chi_1A^*})=\frac{1}{d}\Tr(S^2|_{A^*})$.
\end{lemma}
\begin{proof}
Follows from Lemma~\ref{1112}.
\end{proof}

Let $V_1,\dots, V_n$ be irreducible $A$-modules with characters
$\chi_1,\dots,\chi_n$.
Recall that if $\rho = \sum_j\, \FPdim(V_j) \chi_j$ then 
$\FPdim(A)^{-1}\rho$ is a minimal idempotent
in $R(A)$.

\begin{theorem}[The second trace formula]
\label{the 2nd trace}
We have
\begin{equation}
\label{2nd trace eqn} 
\dim(A) = \frac{\Tr(S^2|_{\rho A^*})}{d} \FPdim(A).
\end{equation}
\end{theorem}
\begin{proof}
Note that the canonical left integral $\ell$ defined by \eqref{ell} 
belongs to the matrix block
corresponding to the trivial representation of $A$ and
$\frac{d}{\Tr(S^2|_{A^*})}\ell$ is a primitive idempotent.
Hence,
$\la \rho,\, \ell\ra = \la \chi_1,\,\ell\ra$. Using 
Lemma~\ref{chi0} we obtain
\begin{eqnarray*}
\Tr(S^2|_{{\rho} A^*})
&=&  \mbox{FPdim}(A)^{-1}  \Tr(\rho\circ S^2|_{A^*}) \\
&=&  \mbox{FPdim}(A)^{-1} \la \rho,\, \ell\ra \\
&=&  \mbox{FPdim}(A)^{-1} \Tr(S^2|{\chi_1 A^*})\\
&=&  \mbox{FPdim}(A)^{-1} d^{-1} \Tr(S^2|_{A^*}),
\end{eqnarray*}
whence the result follows since $\dim(A) =  d^{-2} \Tr(S^2|_{A^*})$.
\end{proof}

\begin{remark}
It was shown in \cite{ENO} that the ratio $\dim(A)/\FPdim(A)$
is an algebraic integer.
\end{remark}

\begin{corollary}
\label{first summand}
We have $\Tr(S^2|_{\rho A^*}) \leq d$.
\end{corollary}
\begin{proof}
It follows from \cite{ENO} that the categorical dimension of $A$
is not bigger than its Frobenius-Perron dimension, 
which implies the statement.
\end{proof}

\begin{corollary}
\label{ps.u = positivity}
A finite-dimensional weak Hopf algebra $A$ is pseudo-unitary 
if and only if the eigenvalues of $S^2$ are positive.
\end{corollary}
\begin{proof}
It was shown in Corollary~\ref{G=F} that for a pseudo-unitary
$A$ the square of the antipode is the conjugation by $wS(w)^{-1}$,
where $w\in Z(A_s)$ has positive eigenvalues.  Since $w$ and $S(w)$ commute
it follows that the eigenvalues of $S^2$ are positive.

To prove the converse, observe that the action of $\rho$ 
on the trivial $A^*$-module $V_1$ is not zero because $\la \rho,\, 1\ra =\mu(A^*)\neq 0$,
where $\mu(A)$ is the index of $A^*$ from Definition~\ref{index}.
Thus $\rho$ gives rise to  a non-zero element in $\End_\mathbb{C}(V_1)$.
Therefore, $\Tr(S^2|_{\rho A^*})\geq \Tr(S^2|_{\rho \End_\mathbb{C}(V_1)}) = d$
and $\dim(A) = \FPdim(A)$ by Corollary~\ref{first summand}.
\end{proof}

\begin{corollary}
The dual of a pseudo-unitary weak Hopf algebra
is pseudo-unitary.
Weak Hopf subalgebras and quotients  of a pseudo-unitary weak Hopf algebra
are pseudo-unitary. 
\end{corollary}

\end{section}


\begin{section}
{The structure of module algebras over weak Hopf algebras}

\subsection{$A$-module ideals and stability of the Jacobson radical}

Let $A$ be a finite-dimensional 
weak Hopf algebra and let $M$ be a finite-dimensional  left $A$-module algebra.

The following Lemma extends \cite[2.2(a)]{Li}.

\begin{lemma}
\label{A-stable}
If $I\subset M$ is a two-sided ideal then 
\begin{equation}
A\cdot I = \left\{ \sum_i\,h_i\cdot x_i\mid h_i\in A,\, x_i\in I\right\}
\end{equation}
is an $A$-stable ideal of $M$.
\end{lemma}
\begin{proof}
It is clear that $A\cdot I$ is an $A$-stable subspace of $M$. 
For all $m\in M,\, h\in A,\, x\in I$ we compute
\begin{eqnarray*}
m(h\cdot x)
&=& (1\1\cdot m)(1\2h\cdot x) = (h\2 S^{-1}(h\1) \cdot m)(h\3 \cdot x) \\
&=& h\2 \cdot ((S^{-1}(h\1)\cdot m)x) \in A\cdot I; \\
(h\cdot x)a
&=& (1\1h\cdot x)(1\2 \cdot a) = (h\1\cdot x)(h\2 S(h\3)\cdot a)\\
&=& h\1 \cdot (x(S(h\2)\cdot a)) \in A\cdot I,
\end{eqnarray*}
and hence, $A\cdot I$ is an ideal.
\end{proof}

\begin{lemma}
\label{LRT}
Let $L$, $R$, and $T$  be linear endomorphisms of $M$
defined by
\begin{equation}
L(m)x = mx,\quad R(m)x = xm,\quad T(h)x= h\cdot x
\mbox{ for all } m,x\in M,\, h\in A.
\end{equation}
Then 
\begin{enumerate}
\item[(i)] for all $y\in A_s$ and $z\in A_t$ we have
$T(y) = R(y\cdot 1)$ and  $T(z) = L(z\cdot 1)$,
\item[(ii)]
for all $h\in A$ and $m\in M$ we have 
$L(h \cdot m) = T(h\1)\circ L(m) \circ T(S(h\2)$.
\end{enumerate}
\end{lemma}
\begin{proof}
We compute 
\begin{equation*}
T(y)m = y\cdot (m1) = (1\1\cdot m)(1\2y \cdot 1) = m (y\cdot 1),
\qquad y\in A_s,\, m\in M,
\end{equation*}
and similarly for $T(z)$, which proves (i). Next, for $h\in A$
and $m,x\in M$ we have 
we compute
\begin{eqnarray*}
T(h\1)\circ L(m) \circ T(S(h\2))x
&=& h\1 \cdot ( m (S(h\2)\cdot x)) \\
&=& (h\1 \cdot m) (h\2 S(h\3) \cdot x) \\
&=& (1\1h \cdot m)(1\2 \cdot x) \\
&=& (h \cdot m)x = L(h \cdot m)x,
\end{eqnarray*}
which proves (ii).
\end{proof}

In \cite{Li} V.~Linchenko proved that the Jacobson radical 
of a module algebra over a semisimple Hopf algebra $A$ is an $A$-stable
ideal. Below we extend this result to pseudo-unitary weak Hopf algebras,
cf.\ Definition~\ref{pseudo-unitary WHA}
(note that a semisimple Hopf algebra is pseudo-unitary).

\begin{theorem}
\label{J(M) is stable}
Let $A$ be a  pseudo-unitary weak Hopf algebra  and let $M$
be a finite-dimensional $A$-module algebra. Then the Jacobson radical $J(M)$ of $M$ 
is an $A$-stable ideal of $M$. 
\end{theorem}
\begin{proof}
We will use notation from Lemma~\ref{LRT}.
Let $G=wS(w^{-1}) \in A$ be the element from Corollary~\ref{G=F}
such that $S^2(h)=GhG^{-1}$ for all $h\in A$. Here $w\in Z(A_s)$ is a Frobenius-Perron
element of $A_s$. In particular, all the eigenvalues of $A_s$ are positive. 
Let $\Tr$ denote the trace of a linear endomorphism of $M$.
We compute
\begin{eqnarray*}
\Tr(T(G)\circ L(h\cdot m))
&=& \Tr(T(S((h\2)G h\1)\circ L(m))) \\
&=& \Tr(T(G S^{-1}(\eps_s(h)))\circ L(m)) \\
&=& \Tr(T(G) \circ L((\eps_s(h)\cdot 1)m))\\
&=& \Tr(L((S(w)\cdot 1)(\eps_s(h)\cdot 1)m)\circ R(w^{-1}\cdot 1)),
\end{eqnarray*}
for all $h\in A,\, m, \in M$. Now if $m\in J(M)$, then
$L((S(w)\cdot 1)(\eps_s(h)\cdot 1)m)$ is a nilpotent endomorphism
of $M$ that commutes with $R(w^{-1}\cdot 1)$, therefore,
\begin{equation}
\Tr(T(G)\circ L(h\cdot m)) =0, \qquad \mbox{ for all } h\in A,\, m\in J(M).
\end{equation}
Since $T(G) = L(w)\circ R(w^{-1})$  and $A\cdot J(M)$
is an  $A$-stable ideal by Lemma~\ref{A-stable} it follows that 
\begin{equation}
\Tr(L(h\cdot m)\circ R(w^{-1})) =0, \qquad \mbox{ for all } h\in A,\, m\in J(M).
\end{equation}
Consequently,
\begin{equation}
\Tr(L(x)^n\circ R(w^{-1})) =0, \qquad \mbox{ for all } x\in A\cdot J(M) \mbox{ and } 
n=1,\,2,\dots.
\end{equation}
Since  $R(w^{-1}))$ is a diagonalizable linear operator with positive eigenvalues
it follows by a standard linear algebra argument that $L(x)$ is a nilpotent
operator for all $x\in A\cdot J(M)$. Hence, any such $x$
is a nilpotent element of $M$ and, therefore,  $A\cdot J(M)$ is a nil ideal of $M$.
It is well known that every nil ideal is contained in the Jacobson radical,
whence the result follows.
\end{proof}

\subsection{Frobenius-Perron dimensions of module algebras}
\label{zhus}

A classical result in group theory states that
if $G$ is a finite group then the cardinality of any transitive $G$-set
divides $|G|$. This was extended to semisimple Hopf algebras by Y.~Zhu
who showed in \cite{Z2} that if $M$ is an indecomposable semisimple module
algebra over a semisimple Hopf algebra $A$ and $M_1,\, M_2$ are irreducible
$M$-modules then $\dim_\mathbb{C}(M)$ divides 
$\dim_\mathbb{C}(M_1)\dim_\mathbb{C}(M_2)\dim_\mathbb{C}(H)$.

In this section we extend Zhu's result to semisimple weak Hopf algebras
and their (co)module algebras. In this case the vector space 
dimensions should be replaced by Frobenius-Perron dimension.

Let $A$ be a semisimple weak Hopf algebra and let 
$M$ be an indecomposable finite-dimensional semisimple left $A$-comodule algebra. 
Recall from Remark~\ref{mod algebra} that in this case
$M$ is  a right $A^*$-module algebra via
\begin{equation}
\label{mod comod}
m\cdot \phi := \la m\I,\, \phi\ra m\II,\qquad m\in M,\,\phi\in A^*,
\end{equation}
where $\delta(m) = m\I \otimes m\II$ is the coaction of $A$.

The Grothendieck group
$K_0(M)$ is a left $K_0(A)$-module via
\begin{equation}
\la \chi\xi,\, m\ra := \la \chi\otimes\xi,\, \delta(m)\ra  = \xi(m\cdot \chi),
\end{equation}
for all $\chi\in K_0(A),\ \xi\in K_0(M)$, and $m\in M$. 

Let $M_1,\dots,M_t$ be irreducible $M$-modules and
let $\xi_1,\dots,\xi_t$ be their characters.

Recall  from Section~\ref{K0 and FP}
that there exists a unique up to a scalar multiple Frobenius-Perron character  
$\xi\in R(M) = K_0(M)\otimes_\mathbb{Z} \mathbb{C}$ of $M$ defined by the property
\begin{equation}
\label{FP M}
\Tr_V\xi = \FPdim(V)\xi,
\end{equation}  
for all finite-dimensional $A$-modules $V$. The coefficients of $\xi$ in the
basis $\{\xi_1,\dots,\xi_t\}$ are Frobenius-Perron dimensions of irreducible
$M$-modules.

The space of Frobenius-Perron characters
of $M$ is $1$-dimensional and is equal to
$\rho R(M)$, where $\rho$ is the Frobenius-Perron
character of $A$. Let $w_{A^*}\in  k\,\eps_s(\rho)$ be any Frobenius-Perron 
element of $A_s^*$ from Definition~\ref{FPdim in As}.
The next Proposition expresses a Frobenius-Perron character of $M$ 
in terms of this element.

\begin{proposition}
\label{xiw}
The element 
\begin{equation}
\xi_f =\sum_{k=1}^t\, \xi_k(1\cdot w_{A^*})\,\xi_k
\end{equation}
is a Frobenius-Perron character of $M$.
\end{proposition}
\begin{proof}
Choose a normalization $w_{A^*} =\eps_s(\rho)$. 
Note that the map 
\begin{equation}
\Phi: K_0(A)\to \mathbb{R} : \xi \mapsto  \xi(1\cdot w_{A^*})
\end{equation}
is a $K_0(A)$-module homomorphism, where $\mathbb{R}$ is a $K_0(A)$-module
via $\FPdim$. Indeed, we have
\begin{eqnarray*}
\Phi(\Tr_V\xi)
&=& \xi(1\cdot w_{A^*}\Tr_V)\\
&=& \xi(1\cdot \rho\Tr_V)\\
&=& \FPdim(V)\xi(1\cdot \rho) = \FPdim(V)\Phi(\xi).
\end{eqnarray*}
Since $\Phi(\xi_k)>0$ for all $k=1,\dots,t$ it follows from Lemma~\ref{unique homo}
that  $\sum_k\,\Phi(\xi_k)\xi_k$ is a Frobenius-Perron character of $M$.
\end{proof}

It follows that we can choose the Frobenius-Perron dimensions of irreducible
$M$-modules to be
\begin{equation}
\FPdim(M_k) = \xi_k(1\cdot w_A^*),\qquad k=1,\dots,t.
\end{equation}

For any $M$-module $U$ with the
character $\xi= \sum_{k=1}^t\, a_k \xi_k \in K_0(M),\, a_k\in \mathbb{Z}_{\geq 0}$
its Frobenius-Perron dimension is
\begin{equation}
\FPdim(U) := \sum_{k=1}^t\, a_k \FPdim(M_k).
\end{equation}
Define the {\em Frobenius-Perron dimension} of $M$ by
\begin{equation}
\FPdim(M) := \sum_{k=1}^t\, \FPdim(M_k)^2.
\end{equation}
Of course, the numbers $\FPdim(M_k)$, $\FPdim(M)$ depend (up to a scalar) 
on our choice, but the ratios 
${\FPdim(M_i)\FPdim(M_k)}/{\FPdim(M)},\, i,k=1,\dots t$ do not.

The following Theorem extends the main result of \cite{Z2}.

\begin{theorem}
\label{orbits}
Let $A$ be a semisimple weak Hopf algebra and $M$ be an indecomposable finite
dimensional semisimple left $A$-comodule algebra. Let $M_1,\dots, M_t$
be irreducible $M$-modules. Then the numbers
\begin{equation}
\frac{\FPdim(A)\FPdim(M_i)\FPdim(M_k)}{\FPdim(M)},\qquad i,k=1,\dots, t,
\end{equation}
are algebraic integers.
\end{theorem}
\begin{proof}
Let $\xi_k$ be the  character
of $M_k, k=1,\dots, t$ and let $e_k$ be the corresponding
primitive central idempotent of $M$. Let 
$\xi_f=\sum_{k=1}^t\,\FPdim(M_k)\xi_k$
be a Frobenius-Perron character of $M$ and
let $\rho$ be the Frobenius-Perron character of $A$. Then
\begin{eqnarray*}
\rho\xi_k 
&=& \sum_j\, \FPdim(V_j) \Tr_{V_j}\xi_k \\
&=& \sum_{ji}\, N_{jk}^i\FPdim(V_j) \xi_i,
\end{eqnarray*}
where $N_{jk}^l$ are some non-negative integers. Hence, $\rho\xi_k$
is a linear combination of $\xi_l,\, l=1\dots, t$ with algebraic integer
coefficients. 
On the other hand, since both $\rho\xi_k$ and $\rho\xi_f$ are Frobenius-Perron
characters of $M$, we have
\begin{eqnarray*}
\rho\xi_k 
&=& \frac{\FPdim(M_k)}{\FPdim(M)}\,\rho\xi_f \\
&=&  \frac{\FPdim(A)\FPdim(M_k)}{\FPdim(M)}\, \xi_f\\
&=& \sum_i\, \frac{\FPdim(A)\FPdim(M_k)\FPdim(M_i)}{\FPdim(M)}\, \xi_i.
\end{eqnarray*} 
Comparing  the coefficients of $\xi_i$ in the two above formulas for $\rho\xi_k$
we get the result.
\end{proof}
\end{section}


\bibliographystyle{ams-alpha}

\begin{thebibliography}{A}  

\bibitem [A]{A} N.~Andruskiewitsch,
\textit{About finite-dimensional Hopf algebras},
in ``Quantum symmetries in theoretical physics and mathematics,''
Contemp.\ Math., \textbf{294}, AMS, Providence (2002), 1-57.

\bibitem [BK]{BK} B.~Bakalov, A.~Kirillov Jr.,
\textit{Lectures on tensor categories and modular functors}, 
AMS, Providence, (2001).

\bibitem [BNSz]{BNSz}  G.~B\"ohm, F.~Nill, and K.~Szlach\'anyi,
\textit{Weak Hopf algebras I. Integral theory and $C^*$-structure},
J. Algebra, \textbf{221} (1999), 385--438. 

\bibitem [BSz1]{BSz1} G.~B\"ohm, K.~Szlach\'anyi,   
\textit{A coassociative $C^*$-quantum group with nonintegral dimensions},
Lett. in Math. Phys, \textbf{35}  (1996), 437--456.

\bibitem [BSz2]{BSz2}  G.~B\"ohm, K.~Szlach\'anyi,
\textit{Weak Hopf algebras II. Representation theory, dimensions
and the Markov trace},
J. Algebra, \textbf{233} (2000), 156-212.

\bibitem [D]{D} V.G.~Drinfeld,
\textit{On almost cocommutative Hopf algebras},
Leningrad Math.\ J.\ \textbf{1} (1990), 321-342.

\bibitem[Ga]{Ga} F.~R.~Gantmacher, \textit{The theory of matrices}, AMS Chelsea
Publishing, Providence, (1998).

\bibitem [GHJ]{GHJ} F.~Goodman, P.~de la Harpe, and V.F.R.~Jones,
\textit{Coxeter graphs and towers of algebras},  M.S.R.I. Publ.\  
\textbf{14}, Springer, Heidelberg, (1989).

\bibitem  [EN]{EN} P.~Etingof, D.~Nikshych,
\textit{Dynamical quantum groups at roots of $1$}, 
Duke Math.\ J., \textbf{108} (2001), 135-168.

\bibitem [ENO]{ENO} P.~Etingof, D.~Nikshych, V.~Ostrik,
\textit{On fusion categories}, \texttt{math.QA/0203060} (2002).

\bibitem [K]{K} G.~Kac,
\textit{Certain arithmetic properties of ring groups},
Functional Anal.\ Appl.\ \textbf{6} (1972), 158-160.

\bibitem [KN]{KN} L.~Kadison,  D.~Nikshych,
\textit{Frobenius extensions and weak Hopf algebras},
J.\ Algebra, \textbf{244} (2001), 312-342.

\bibitem [LR1]{LR1} R.~Larson, D.~Radford,
\textit{Finite-dimensional cosemisimple Hopf algebras in characteristic $0$ 
are semisimple}, J. Algebra \textbf{117} (1988), 267--289.

\bibitem [LR2]{LR2} R.~Larson, D.~Radford,
\textit{Semisimple cosemisimple Hopf algebras}, 
Amer.\ J.\ Math.\  \textbf{109} (1987), 187--195.

\bibitem [Li]{Li} V.~Linchenko,
\textit{Nilpotent subsets of Hopf module algebras},
to appear in ``Groups, rings, Lie and Hopf algebras'',
Mathematics and its applications, \textbf{555} (2003).

\bibitem [Lo]{Lo} M.~Lorenz, 
\textit{On the class equation for Hopf algebras},
Proc.\ AMS, \textbf{126} (1998), 2841-2844.

\bibitem [Ma]{Ma} S.~MacLane,
\textit{Categories for the working mathematician},
2nd edition, Springer-Verlag (1998). 

\bibitem [Mo]{Mo} S.~Montgomery,
\textit{Hopf algebras and their actions on rings}, CBMS Regional Conference
Series in Mathematics, \textbf{82}, AMS, (1993).

\bibitem [Mo1]{Mo1} S.~Montgomery,
\textit{Classifying finite-dimensional semisimple Hopf algebras},
in ``Trends in the representation theory of finite-dimensional algebras''
(Seattle, 1997), 256-279, Contemp.\ Math.\, \textbf{229}, AMS, (1998).

\bibitem [Mo2]{Mo2} S.~Montgomery,
\textit{Representation theory of semisimple weak Hopf algebras},
Algebra--Representation theory (Constanta 2000), 189-218,
NATO Sci.\ Ser.\ II Math.\ Phys.\ Chem., \textbf{28} (2001).

\bibitem [Mu]{Mu} M.~M\"uger, 
\textit{From subfactors to categories and topology I. Frobenius algebras in 
and Morita equivalence of tensor categories}, 
J. Pure Appl. Algebra \textbf{180} (2003), 81-157.

\bibitem [N1]{N1} D.~Nikshych,
\textit{Duality for actions of weak Kac algebras and crossed product
inclusions of II$_1$ factors}
J.\ Operator Theory \textbf{46}, 635-655 (2001).

\bibitem [N2]{N2} D.~Nikshych,
\textit{On the structure of weak Hopf algebras},
Adv.\ Math.\ \textbf{170}, 257-286 (2002).

\bibitem [NTV]{NTV} D.~Nikshych, V.~Turaev, and L.~Vainerman,
\textit{Quantum groupoids and invariants of knots and 3-manifolds},
Topology and Applications,  \textbf{127} (2003), 91-123.

\bibitem  [NV1]{NV1} D.~Nikshych, L.~Vainerman,
\textit{A characterization of depth $2$ subfactors of II${}_1$  factors},
J.\ Func.\ Analysis, \textbf{171} (2000), 278-307.

\bibitem [NV2]{NV2} D.~Nikshych, L.~Vainerman,
\textit{Finite quantum groupoids and their applications},
in ``New Directions in Hopf Algebras,''
MSRI Publications, \textbf{43} (2002), 211-262.

\bibitem [O]{O}  V.~Ostrik,
\textit{Module categories, weak Hopf algebras and modular invariants},
Transform. Groups, \textbf{8} (2003), 177-206.

\bibitem [R]{R} D.~Radford,
\textit{The order of the antipode of a finite-dimensional Hopf algebra 
is finite}, Amer.\ J.\ Math., \textbf{98} (1976), 333--355.

\bibitem [S]{S} D.~Stefan,
\textit{The set of types of $n$-dimensional semisimple and cosemisimple 
Hopf algebras is finite}, J. Algebra \textbf{193} (1997), no. 2, 571--580.

\bibitem [V]{V} P.~Vecsernyes,
\textit{Larson-Sweedler theorem, grouplike elements, and invertible
modules in weak Hopf algebras}, preprint (2001), \texttt{math.QA/0111045}.

\bibitem [Z1]{Z1} Y.~Zhu,
\textit{Hopf algebras of prime dimension},
Internat.\ Math.\ Res. Notices (1994), no. 1, 158-160.

\bibitem [Z2]{Z2} Y.~Zhu,
\textit{The dimension of irreducible modules for transitive
Hopf algebras}, Comm.\ Alg.,\ \textbf{29} (2001), 2877-2886.

\end{thebibliography}
  
\end{document}